\documentclass[a4paper,fullpage]{article}
\usepackage{amssymb}
\usepackage{indentfirst}
\usepackage{amsmath,amssymb,amsthm,color}
\usepackage{latexsym}
\usepackage{sidecap}
\usepackage{hyperref}
\usepackage[applemac]{inputenc}
\usepackage{marginnote}
\usepackage[T1]{fontenc}

\newcommand{\R}{\mathbb{R}}

\renewcommand{\theta}{\vartheta}
\newcommand{\al}{\alpha}

\providecommand{\keywords}[1]{{\small{\sc Keywords: }#1}}
\providecommand{\msc}[1]{{\small{\sc MSC: }#1}}

\title{A countable number of homoclinic loops as an $\omega$-limit set for a planar vector field \date{}} 
\author{Francesco Spadaro\footnote{Bachelor student of Mathematics at {\sc Dipartimento di Matematica ``G. Castelnuovo'', Sapienza Università di Roma,} P.le Aldo Moro 5, 00185 Roma, Italy. {\it email:}{\tt fspadaro.cisco@gmail.com}}}  

\numberwithin{equation}{section}
\newtheorem{theorem}{Theorem}[section]

\theoremstyle{definition}
\newtheorem{defin}{Definition}[section]
\newtheorem{remark}{Remark}[section]

\begin{document}
\maketitle

\begin{abstract}
We construct an example of a smooth, non analytic planar vector field with an $\omega$-limit set consisting of a countable number of homoclinic loops.
\end{abstract}
\bigskip
\\ \\
\msc{34C05, 34C35, 34C37, 34C40, 58F13}
\section{Introduction}
In the last two centuries the behavior of planar dynamical systems has been deeply investigated and understood. 

The Poincaré-Bendixson theorem plays a fundamental role in describing the long-time behavior of solutions of planar dynamical systems.

This theorem completely characterizes the asymptotic behavior of solutions of planar dynamical systems having a finite number of equilibrium points. Under this last assumption the heteroclinic orbits need to be in a finite number too, while the homoclinic orbits can be even infinite. In this article we recall under which hypothesis this may occur and we explicitly build a system with such a behavior.

I thank my advisor Prof. M. A. Pozio for her help and Prof. P. Buttà e Prof. P. Negrini of Sapienza, Università di Roma, for their suggestions.

\section{Preliminaries}
Let us introduce the notation used in this article.

Let $\phi:S\times\R\to S$, $\phi:(x,t)\mapsto\phi^{t}(x)$ be the flux for the Cauchy problem
\begin{equation}\label{cauchy}
\begin{cases}
\dot x =f(x),
\\ x(t_{0})=x_{0}\in S.
\end{cases}
\end{equation}
where $S\subseteq\R^{2}$ is an open set and $f:S\to\R^{2}$ is a $C^{1}(S)$ vector field. Moreover $
\gamma^{+}(x_{0})=\bigcup_{t\geq0}\phi^{t}(x_{0})$ is called \emph{the positive semi-orbit} ($\gamma^{-}(x_{0})=\bigcup_{t\leq0}\phi^{t}(x_{0})$ \emph{the negative semi-orbit}) of $x_{0}$, and
$\gamma=\gamma(x_{0})=\gamma^{+}(x_{0})\cup\gamma^{-}(x_{0})$ is called \emph{the orbit} of $x_{0}$.

\begin{defin}
Given an orbit $\gamma\subset S$, the set $\omega(\gamma)=\bigcap_{x_{0}\in\gamma}\overline{\gamma^{+}(x_{0})}$ is called \emph{the $\omega$-limit set} of $\gamma$, similarly $\al(\gamma)=\bigcap_{x_{0}\in\gamma}\overline{\gamma^{-}(x_{0})}$ is called \emph{the $\al$-limit set} of $\gamma$.
\end{defin}
The $\omega$-limit set (respectively the $\al$-limit set) does not depend on a particular point of the orbit $\gamma$ and every point $p\in S$ belongs to exactly one orbit $\gamma$. Thus we will use both the notation $\omega(\gamma)$ and $\omega(p)$ (respectively $\al(\gamma)$ and $\al(p)$).
\begin{defin}
A set $D\subseteq S$ is  a {\it positive} (respectively {\it negative}) {\it invariant set} for the flux $\phi$ if $\phi^{t}(D)\subseteq D$ for all $t\geq0$ (respectively $t\leq0$). If a set is both positive and negative invariant then it is an {\it invariant set}.
\end{defin}
\begin{defin}
An orbit $\gamma$ that connects two different equilibrium points $x_1,\,x_2$, i.e. such that $\alpha(\gamma)= x_1$ and $\omega(\gamma)= x_2$ is called a {\it heteroclinic orbit}. A non trivial orbit such that $x_1= x_2$ is called a {\it homoclinic orbit} (or loop).
\end{defin}

We now recall the Poincaré-Bendixson theorem.
\begin{theorem}[Poincaré-Bendixson]\label{pb}
Let $D\subseteq S$ be a positively invariant, simply connected and compact region containing a finite number of equilibrium points. Then one of the following holds:
\begin{enumerate}
\item $\omega(p)$ is an equilibrium point;
\item $\omega(p)$ is a closed orbit;
\item $\omega(p)$ is the union of equilibrium points and a finite number of heteroclinic orbits and/or a possibly countable number of homoclinic orbits.
\end{enumerate}
\proof
For a proof of this theorem see \cite[page 18]{palis}. Concerning ``{\it 3.}'', the authors of \cite{palis} prove that given two equilibria, only one heteroclitic orbit may exist, having the first equilibrium as $\omega$-limit and the second as $\alpha$-limit. Hence the number of heteroclinic orbits is finite \cite[page 17]{palis}. But the authors of \cite{palis} do not prove that the $\omega$-limit set is the union of an at most countable number of homoclinic orbits, then we prove it here.

Since we deal with autonomous systems, by uniqueness, orbits cannot intersect. Then if a homoclinic loop is approached from outside, no other loops inside it can be in the $\omega$-limit set. Moreover the orbit starting at $p$ can approach at most one homoclinic loop of the set $\omega(p)$ from inside. Let us prove that it can approach at most a countable number of homoclinic loops from outside. Indeed $D$ is a compact region, thus it has a finite Lebesgue measure, $\mu(D)$. By the Jordan curve theorem, any homoclinic orbit $\gamma$ encloses a bounded region $P(\gamma)$ with a strictly positive Lebesgue measure $\mu(P)$. Moreover the regions corresponding to different homoclinic loops may intersect at the equilibrium only. For every integer $n \geq1$, since $\mu(D)$ is finite, there can be no more than $n$ homoclinic orbits $\gamma$ such that
\begin{equation}\label{measure}
\frac{\mu(D)}{n+1}<\mu(P(\gamma))\leq\frac{\mu(D)}{n}
\end{equation}
Then the number of homoclinic orbits in the $\omega$-limit set is at most countable.

\endproof
\end{theorem}

\begin{remark}
If $D$ contains a single equilibrium point $z$ which is hyperbolic for the field $f$, i.e. the two eigenvalues of $\nabla f(z)$ have non zero real part, by the stable manifold theorem there exist at most two homoclinic loops that close in $z$. But, if $z$ is not a hyperbolic equilibrium point, it might happen to have a larger number of homoclinic loops. However, as far as the vector field is analytic, the number of homoclinic orbits has to be finite \cite[p. 30]{anosov}.
\end{remark}
\begin{remark}
It should be pointed out that also a continuum of homoclinic orbits relative to the same equilibrium may exist, as it is the case in the example we construct below. However the $\omega$-limit set may be the union of an at most countable number of them, as we proved above.
\end{remark}

In \cite{anosov} it is observed that there may exist a planar dynamical system with $\omega$-limit sets consisting of an infinite number of homoclinic loops, though in this case the vector field cannot be analytic. Indeed, in \cite[p. 186]{sansone} an example of such a field is proposed, though it is only a once differentiable field.
Instead, we construct here an explicit example of an infinitely differentiable field, so to underline that the pathology is in the non-analyticity of the field only.

\section{Constructing the vector field}
For constructing the aimed vector field we switch to polar coordinates centered at the equilibrium point, say $(0,0)$. Then it is much easier to identify the attractive manifold with infinite homoclinic loops that we are looking for.

We look for two functions $\psi:\R^{+}\to\R^{+}$ and $F:\R^{+}\times\R\to\R$, $F\in C^{1}(\R^{+}\times\R)$, $2\pi$-periodic in the second variable (here $\R^{+}=[0,+\infty)$). Then we consider the dynamical system
\begin{equation}\label{pvf}
\begin{cases}
\dot r=-\psi(r)\left(F(r,\theta)+\frac{\partial F}{\partial\theta}(r,\theta)\right),
\\ \dot \theta=\,\,\,\,\psi(r)\frac{\partial F}{\partial r}(r,\theta).
\end{cases}
\end{equation}

Along the solutions of system \eqref{pvf} we have
\begin{equation}\label{dotf}
\dot F=\frac{d}{dt}F(r(t),\theta(t))=-\psi F\partial_{r}F
\end{equation}
Therefore, if the Cauchy problem associated to \eqref{pvf} has an initial value $(r_{0},\theta_{0})$ such that $F(r_{0},\theta_{0})=0$, then $F(r(t),\theta(t))=0$ along the solution, or equivalently the solution for this Cauchy problem lies on the graph of $F(r,\theta)=0$.
Thus we look for a function $F:\R^{+}\times\R\to\R$ such that the graph of $F(r,\theta)=0$ has an infinite number of loops.
A suitable function is
\begin{equation}\label{effe}
F(r,\theta)=r-e^{-\frac{1}{(1+\cos\theta)^{2}}}\sin^{2}(\tan(\theta/2)).
\end{equation}
The function is not well-defined for $\theta=\pi+2k\pi, k\in\mathbb{Z}$ but it can be extended by continuity there. From now on $F$ will be the extended function. The equation $F(r,\theta)=0$ reads
\begin{equation}
r=e^{-\frac{1}{(1+\cos\theta)^{2}}}\sin^{2}(\tan(\theta/2)),
\end{equation}
since the right hand side assumes the zero value an infinite number of times, the graph of $F(r,\theta)=0$ gives the {\it edge of a rose} with an infinite number of petals. Let us define it as
\begin{equation}
\Omega:=\{(x,y)\in\R^{2}:F(r(x,y),\theta(x,y))=0\}.
\end{equation}
We choose for the vector field in \eqref{pvf} the function
\begin{equation}
\psi(r)=
\begin{cases}
e^{-\frac{1}{r^{2}}} \qquad\qquad r\neq0
\\ \,\,\,\,0 \qquad\qquad\quad r=0
\end{cases}
\end{equation}
then we explicitly write the vector field in polar coordinates:
\begin{equation}\label{sishomo}
\begin{cases}
\dot r=-e^{-\frac{1}{r^{2}}}H(r,\theta),\\
\dot \theta=e^{-\frac{1}{r^{2}}}.
\end{cases}
\end{equation}
where $H(r,\theta)=r\geq0$ if $\theta=\pi+2k\pi, k\in\mathbb{Z}$ otherwise 
$$H(r,\theta)=r-e^{-\frac{1}{(1+\cos\theta)^{2}}}\left(\sin^{2}(\tan(\theta/2))-\frac{2\sin^{2}(\tan(\theta/2))\sin\theta}{(1+\cos\theta)^{3}}+\frac{\sin(2\tan(\theta/2))}{1+\cos\theta}\right),$$
and in cartesian coordinates:
\begin{equation}\label{sishomocart}
\begin{cases}
\dot x= -e^{-\frac{1}{r^{2}}}\left(y+x-\frac{x}{r}G(x,y)\right),
\\
\dot y= -e^{-\frac{1}{r^{2}}}\left(y-x-\frac{y}{r}G(x,y)\right),
\end{cases}
\end{equation}
where,  if $ y\neq0$ or $y=0$ and $x>0$, we set
$$G(x,y)=e^{-(\frac{r}{r+x})^{2}}\left(\sin^{2}\left(\frac{y}{r+x}\right)-\frac{2yr^{2}\sin^{2}\left(\frac{y}{r+x}\right)}{(r+x)^{3}}+\frac{r}{r+x}\sin\left(\frac{2y}{r+x}\right)\right),$$
otherwise we set $G(x,y)=0$.

The vector field in \eqref{sishomocart} may fail to be smooth on the negative $x$ semi-axis. However for any given $x_{0}\leq0$ we have
\begin{equation}
\lim_{(x,y)\to(x_{0},0)}\frac{e^{-\frac{1}{r^{2}}}e^{-(\frac{r}{r+x})^{2}}}{(r+x)^{n}}=\lim_{(x,y)\to(x_{0},0)}\frac{e^{-\frac{1}{r^{2}}}}{r^{n}}\Big(\frac{r}{r+x}\Big)^{n}e^{-(\frac{r}{r+x})^{2}}=0,
\end{equation}
where for $x_{0}=0$ we have used the boundedness of $\Big(\frac{r}{r+x}\Big)^{n}e^{-(\frac{r}{r+x})^{2}}$. Thus, the vector field in \eqref{sishomocart} is $C^{\infty}(\R^{2})$ since it is continuous with all its derivatives.

Let us prove that the infinite petals of the rose $\Omega$ are the $\omega$-limit set of some orbits of system \eqref{sishomocart}. Indeed, the following holds for system \eqref{sishomocart} having a $C^{\infty}(\R^{2})$ vector field.

\begin{theorem}
For any $(x_{0},y_{0})$, let $\gamma$ be the orbit for system \eqref{sishomocart} starting at $(x_{0},y_{0})$. Then
\begin{enumerate}
\item If $F(r(x_{0},y_{0}),\theta(x_{0},y_{0}))\leq0$, $\gamma$ is a homoclinic loop and $\omega(\gamma)$, $\al(\gamma)$ coincide with the origin;
\item if $F(r(x_{0},y_{0}),\theta(x_{0},y_{0}))>0$, $\gamma$ has $\Omega$, the {\emph edge of a rose with a countable number of petals}, as its $\omega$-limit set.
\end{enumerate}
\begin{remark}
The result in part ``{\it 1.}'' of the last theorem proves that each petal inside $\Omega$ is filled with a continuum of homoclinic orbits.
\end{remark}
\proof
First of all, let 
\begin{equation}
r_0= r(x_0,y_0)= \sqrt{x_0^2+y_0^2}, \quad\quad \theta_0 = \arg(x_0,y_0),\quad\quad F_{0}=F(r_{0},\theta_{0}).
\end{equation}
By \eqref{dotf} the case $F(r_{0},\theta_{0})=0$ follows easily and the orbit is a homoclinic loop lying on $\Omega$.

If $F(r_{0},\theta_{0})<0$, the orbit starts inside a petal and it cannot leave the petal since orbits cannot cross. System \eqref{sishomo}, or equivalently \eqref{sishomocart}, has only one equilibrium point, the origin, which is the $\omega$-limit set if we prove that ``{\it 2.}'' and ``{\it 3.}'' of Theorem \ref{pb} cannot occur. Indeed an eventually closed orbit needs to enclose at least an equilibrium point \cite[p. 251]{hirsch}, hence this cannot be the case since no other equilibria exist inside the petals. Moreover if the orbit approaches the petal from inside, as its $\omega$-limit set, the $\alpha$-limit set would be an equilibrium inside the petal, but this is not possible. Hence the conclusion. The same arguments work for the $\alpha$-limit set.

We prove the second part of the theorem showing that the origin cannot be the $\omega$-limit set of points outside the petals for which $F_{0}>0$. 
Without regularization, say for $\psi=1$, \eqref{dotf} reads
\begin{equation}
\frac{d}{dt}F(r(t),\theta(t))=-F \quad\Rightarrow\quad F(r(t),\theta(t))=F_{0}e^{-t},
\end{equation}
where $F_{0}>0$. Then, $F$ is positive for all $t\in[0,\infty)$ and vanishes as $t$ goes to infinity, therefore orbits need to tend towards $F(r,\theta)=0$, i.e. $\Omega$.
Furthermore, for $\psi=1$, we have $\theta(t)=t$, hence the orbits twist around the rose an infinite number of times, getting closer and closer to it. The function $\psi$, that we used to regularize the field, changes the amplitude of the field but not its direction, except for the origin. Then the orbits do not change, so the rose is, indeed, the $\omega$-limit set.
\endproof
\end{theorem}

\end{document}